\documentclass[reqno]{amsart}
\usepackage{pictex}
\usepackage{mathrsfs}
\usepackage{hyperref}
\usepackage{color}
\begin{document}

%
%

\def\labelenumi{(\theenumi)}

\newtheorem{thm}{Theorem}[section]
\newtheorem{lem}[thm]{Lemma}
\newtheorem{conj}[thm]{Conjecture}
\newtheorem{cor}[thm]{Corollary}
\newtheorem{add}[thm]{Addendum}
\newtheorem{prop}[thm]{Proposition}
\theoremstyle{definition}
\newtheorem{defn}[thm]{Definition}
\theoremstyle{remark}
\newtheorem{rmk}[thm]{Remark}
\newtheorem{example}[thm]{{\bf Example}}

\newcommand{\OmegaH}{\Omega/\langle H \rangle}
\newcommand{\hatOmegaHstar}{\hat \Omega/\langle H_{\ast}\rangle}
\newcommand{\SurfG}{\Sigma_g}
\newcommand{\TriangG}{T_g}
\newcommand{\TriangGOne}{T_{g,1}}
\newcommand{\ProjG}{\mathcal{P}_g}
\newcommand{\TeichG}{\mathcal{T}_g}
\newcommand{\CirclePackGTau}{\mathsf{CPS}_{g,\tau}}
\newcommand{\CrossRatio}{{\bf c}}
\newcommand{\CrossRatioGTau}{\mathcal{C}_{g,\tau}}
\newcommand{\CrossRatioOneTau}{\mathcal{C}_{1,\tau}}
\newcommand{\DeformGTau}{\mathcal{C}_{g,\tau}}
\newcommand{\Forget}{\mathit{forg}}
\newcommand{\Uniform}{\mathit{u}}
\newcommand{\Section}{\mathit{sect}}
\newcommand{\SLTwoC}{\mathrm{SL}(2,{\mathbb C})}
\newcommand{\SLTwoR}{\mathrm{SL}(2,{\mathbb R})}
\newcommand{\SUTwo}{\mathrm{SU}(2)}
\newcommand{\PSLTwoC}{\mathrm{PSL}(2,{\mathbb C})}
\newcommand{\GLTwoZ}{\mathrm{GL}(2,{\mathbb Z})}
\newcommand{\GLTwoC}{\mathrm{GL}(2,{\mathbb C})}
\newcommand{\PSLTwoR}{\mathrm{PSL}(2,{\mathbb R})}
\newcommand{\PGLTwoR}{\mathrm{PGL}(2,{\mathbb R})}
\newcommand{\GLTwoR}{\mathrm{GL}(2,{\mathbb R})}
\newcommand{\PSLTwoZ}{\mathrm{PSL}(2,{\mathbb Z})}
\newcommand{\SLTwoZ}{\mathrm{SL}(2,{\mathbb Z})}
\newcommand{\nnn}{\noindent}
\newcommand{\MCG}{{\mathcal {MCG}}}
\newcommand{\MMap}{{\bf \Phi}_{\mu}}
\newcommand{\HH}{{\mathbb H}^2}
\newcommand{\TT}{{\mathbb T}}
\newcommand{\X}{{\mathcal  X}}
\newcommand{\C}{{\mathscr C}}
\newcommand{\CC}{{\mathbb C}}
\newcommand{\RR}{{\mathbb R}}
\newcommand{\Q}{{\mathbb Q}}
\newcommand{\ZZ}{{\mathbb Z}}
\newcommand{\PL}{{\mathscr {PL}}}
\newcommand{\GP}{{\mathcal {GP}}}
\newcommand{\GT}{{\mathcal {GT}}}
\newcommand{\GQ}{{\mathcal {GQ}}}
\newcommand{\EE}{{{\mathcal E}(\rho)}}
\newcommand{\HHH}{{\mathbb H}^3}
\def\square{\hfill${\vcenter{\vbox{\hrule height.4pt \hbox{\vrule width.4pt
height7pt \kern7pt \vrule width.4pt} \hrule height.4pt}}}$}

\newenvironment{pf}{\noindent {\sl Proof.}\quad}{\square \vskip 12pt}

\title{ A simple proof of the Markoff conjecture for prime powers}
\author{Mong Lung Lang, and Ser Peow Tan}
\address{Department of Mathematics \\ National University of Singapore \\
2 Science Drive 2 \\ Singapore 117543} \email{matlml@nus.edu.sg;
mattansp@nus.edu.sg}
 \vskip 20pt

%
%

\begin{abstract}
We give a simple and independent proof of the result of Jack
Button and Paul Schmutz that the Markoff conjecture on the
uniqueness of the Markoff triples $(a,b,c)$ where $a \le b \le c$
holds whenever $c$ is a prime power.
\end{abstract}

 \maketitle
\vskip 6pt

\section{Introduction}\label{s:intro}

The Markoff conjecture (first conjectured by G. Frobenius in 1913)
states that the set of triples $(a,b,c)$ of positive integer
solutions of the Markoff equation
\begin{equation}\label{eqn:markoff}
    a^2+b^2+c^2=3abc
\end{equation}
is uniquely determined by $c$, if we order the triple so that $a
\le b \le c$. The triples $(a,b,c)$ are called \emph{Markoff
triples} and the values \emph{Markoff numbers}.

\begin{thm}\label{thm:maintheorem}{\rm (}J.Button,
P. Schmutz {\rm )} The Markoff conjecture is true when $c$ is a
prime power.

\end{thm}

This was proven independently  by Schmutz \cite{Schmutz1996} and
Button \cite{Button1998}. We give here a very short and simple
proof of their result which depends only on some elementary
hyperbolic geometry and congruences. The main novelty in our
method is that it combines the geometric and arithmetic
information in a hitherto unexplored way.

The interest in the conjecture lies in the fact that it ties in
several apparently unrelated fields including diophantine
approximations, quadratic forms and hyperbolic geometry, the
reader is referred to \cite{CusickFlahive} for background material
on the number theoretic aspects. The relation with hyperbolic
geometry and the trace field of the modular torus, in particular,
the relation between (\ref{eqn:markoff}) and the Fricke trace
identities appears to be first discovered by H. Cohn in
\cite{Cohn}, see \cite{Series} for an excellent survey of the
subject.

\vskip 10pt

 \noindent {\it Acknowledgements.}
We obtained our proof soon after learning of Button's result in
his Warwick PhD thesis. We had not written it out then as there
was no essentially new results, and it was clear that the method
alone would not resolve the Markoff conjecture in full.
Nonetheless, there seemed to be some interest in the proof,  we
thank in particular Caroline Series and Brian Bowditch for their
interest and encouragement.
\section{Some basic facts}\label{s:basicfacts}
We state here some basic facts which seem to be well-known to the
experts in the field, and also some useful propositions.

\subsection{Markoff triples and the Markoff
tree}\label{ss:markofftriples}
 The Markoff triples $(a,b,c)$ can
be generated from the basic triple $(1,1,1)$ as values attached to
the vertices of an infinite binary tree via the operations
$$(a,b,c) \mapsto (b,c,3bc-a),\quad (a,b,c) \mapsto (a, 3ac-b,
c),\quad (a,b,c) \mapsto (a,b,3ac-b).$$ (It is perhaps visually
more pleasing and systematic to think of the Markoff numbers as
values attached to the complementary regions of the tree, and the
triples as arising from the three regions adjacent to a vertex of
the tree, see for example \cite{Bowditch}). There is a $D_6$
symmetry for this tree of values and the conjecture is equivalent
to saying that the values are unique up to the action of $D_6$.

\subsection{The modular torus and simple closed geodesics}\label{ss:modulartorus}
The modular torus is defined to be the cusped hyperbolic torus
${\mathbb T}=\mathbb H/G$, where $G=\langle A,B \rangle
<\SLTwoZ=\Gamma$,
$$A=\left(%
\begin{array}{cc}
  2 & 1 \\
  1 & 1 \\
\end{array}%
\right), \quad B=\left(%
\begin{array}{cc}
  1 & 1 \\
  1 & 2 \\
\end{array}%
\right)\in \Gamma.$$ Note that $G=[\,\Gamma,\Gamma]$, see for
example \cite{LLT}. Also, it is common to regard $A,B$ as matrices
in $\PSLTwoZ$, we have chosen a lift to $\SLTwoZ$ here for which
the traces of all simple closed geodesics are positive. The
isometry group of ${\mathbb T}$ is also isomorphic to $D_6$,
 and there is a direct correspondence between the Markoff numbers and
the traces/lengths of simple closed geodesics on ${\mathbb T}$,
indeed, the traces are precisely 3 times the Markoff numbers.
Similarly, there is a correspondence between the Markoff triples
and the traces/lengths of triples of simple closed geodesics on
${\mathbb T}$ with pair-wise geometric intersection number one
(see \cite{Cohn}, \cite{Series} or \cite{Bowditch}). The Markoff
conjecture is equivalent to saying that the traces/lengths of the
simple closed geodesics on $\mathbb T$  are distinct, up to the
action of the isometry group $D_6$ on $\mathbb T$.

\subsection{Markoff Matrices}\label{ss:markoffmatrices}
A matrix $M \in G$ is \emph{primitive} if  it corresponds to a
primitive curve on $\mathbb T$, that is, $M \neq N^n$ for some $N
\in G$, $n \neq \pm 1$. A primitive matrix $M \in G$ is called a
\emph{Markoff matrix} if (i) it corresponds to a simple closed
geodesic $\gamma$ on ${\mathbb T}$; and (ii) its fixed axis
corresponds to a lift of $\gamma$ with maximum height in $\HH$.
Note that if $M$ is a markoff matrix, then so is $M^{-1}$.

Let $T_{\alpha}=\left(%
\begin{array}{cc}
  1 & \alpha \\
  0 & 1 \\
\end{array}%
\right)$, and use $T$ for $T_1$. Then
$[A,B^{-1}]=AB^{-1}A^{-1}B=-T_6$. We have the following which
again is a reinterpretation of well-known results on the simple
closed geodesics on $\mathbb T$.

\begin{prop}\label{prop:Markoffmatrices}
If  $M=\left(%
\begin{array}{cc}
  a & b \\
  c & d \\
\end{array}%
\right)$ is a Markoff matrix, then $|c|$ is a Markoff number and
$a+d=3|c|$. Conversely, for any Markoff number $c$, there is a
Markoff matrix $M$ with ${\rm tr}\,M=3c$ and the $(2,1)$ entry of
$M$ equal to $c$. Furthermore, two Markoff matrices $M$ and $N$
correspond to the same simple closed geodesic in ${\mathbb T}$ if
and only if $M^{\pm 1}= T^{3n}NT^{-3n}$ for some $n \in \mathbb
Z$.
\end{prop}

Let $M'$ be obtained from $M$ by interchanging the diagonal
entries. Taking into consideration the $D_6$ symmetry of $\mathbb
T$, we define an equivalence relation $\sim$ on the set of Markoff
matrices by $M \sim N$ if and only if $N^{\pm 1}=T^nMT^{-n}$, or
$N^{\pm 1}=T^nM'T^{-n}$ for some $n \in \mathbb Z$.  Since the
$D_6$ symmetry of $\mathbb T$ is generated by conjugation by
reflection on the hyperbolic line $(0, \infty)$ and conjugation by
$T$, two Markoff matrices are equivalent if and only if they
either correspond to the same simple closed geodesic, or to two
simple closed geodesics on $\mathbb T$ equivalent under the action
of $D_6$. The Markoff conjecture is hence equivalent to the
statement that the  equivalence classes of Markoff matrices are
uniquely determined by the traces, or alternatively, $|c|$, where
$c$ is the $(2,1)$ entry. This is the statement we will prove, in
the case $|c|=p^n$, where $p$ is prime, in the next section. We
first prove some simple technical propositions on Markoff
matrices.

For a Markoff matrix $M$, and each $k \in \mathbb Z$, define
$$M_k=
\left (
\begin{array}{ll}
1 & k/c \\
0 & 1 \\
\end{array}
\right ) \left (
\begin{array}{ll}
a & b \\
c & d \\
\end{array}
\right ) \left (
\begin{array}{ll}
1 & -k/c \\
0 & 1 \\
\end{array}
\right ) = \left (
\begin{array}{ll}
a+k & b+k(d-a-k)/c \\
c & d-k \\
\end{array}
\right )\,.
$$ It is clear that $M_k \in \SLTwoZ $
 if and only if $k(d-a-k)$ is a multiple of $c$.

\begin{prop}\label{prop:cnotmultipleof4}
Let $M=\left(%
\begin{array}{cc}
  a & b \\
  c & d \\
\end{array}%
\right)$  be a Markoff matrix with $c>0$.
 Then $c$ is not a multiple of $4$ and all the odd prime divisors
 of $c$ are of the form $4m+1$ for some $m \in  \Bbb N$.

\end{prop}

\begin{pf}
Since  $1 = ad -bc =3cd-d^2-bc$, we have  $ d^2 \equiv -1$
(mod $c$). This implies that
 $$ X^2 \equiv -1 \,(\mbox{ mod } c )$$
is solvable. It follows that $c$ is not
 a multiple of $4$ and possesses no prime divisor
$p$ such that $p \equiv 3$ (mod 4).

\end{pf}

\medskip

\begin{prop}\label{prop:exactdivisorofc} Let $M=\left(%
\begin{array}{cc}
  a & b \\
  c & d \\
\end{array}%
\right)$  be a Markoff matrix and
 let $k \in \Bbb Z$.
Then $gcd\,(c,k, d-a-k)  = 1$ or $2$. In particular, if $M_k \in
 \SLTwoZ $, then any
exact divisor  $p^m$ of $c$ $(p$ is an odd prime$)$
 is relatively prime to either  $k$ or $d-a-k$.

\end{prop}

\begin{pf}
 Suppose
  gcd$\,(c,k, d-a-k) \ne 1$.  Let $x$ be a divisor
 of $(c,k, d-a-k)$.
 It follows that $x$ is a divisor of $d-a$.
Since $a+d = 3c$ and $x$ divides $c$,
 $x $ is a divisor of $a+d$.
Consequently, $x$ is a divisor of $2a= (a-d)+(a+d)$. Since $x|c$,
$x|2a$ and  $(a,c) = 1$, we conclude
 that $x = 2$.
\end{pf}

\section{Proof of Theorem}\label{s:prrofoftheorem}
Theorem \ref{thm:maintheorem} is now equivalent to the following:
\begin{lem}\label{lem:uniquenessofmarkoffmatrics}
Suppose that $M$ and $N$ are two Markoff matrices with ${\rm
tr}\,M={\rm tr}\,N=3p^n$ where $p$ is prime. Then $M \sim N$.
\end{lem}

\begin{pf}
Let $M=\left(%
\begin{array}{cc}
  a & b \\
  c & d \\
\end{array}%
\right), \quad N=\left(%
\begin{array}{cc}
  a' & b' \\
  c' & d' \\
\end{array}%
\right)$. Taking inverses if necessary, we may assume that
$c=c'>0$, and so $a+d=a'+d'=3c$. Hence $a'=a+k, d'=d-k$ for some
$k \in \mathbb Z$. The key observation here is that since $M$ and
$N$ have the same traces and their fixed axes have the same
height,  $M$ and $N$ are conjugate by a parabolic transformation
fixing $\infty$. A simple computation gives $$N=M_k=T_{k/c}\cdot
M\cdot T_{k/c}^{-1}=\left (
\begin{array}{ll}
a+k & b+k(d-a-k)/c \\
c & d-k \\
\end{array}
\right )\,.$$
%
Since $N \in \SLTwoZ $, $k(d-a-k)$ is a multiple of $c$. By
proposition \ref{prop:cnotmultipleof4}, we may assume that $p$ is
a odd prime, and by proposition \ref{prop:exactdivisorofc},
$c=p^n$ divides $k$ or $(d-a-k)$. In the first case,
$N=T^lMT^{-l}$ for some $l \in \mathbb Z$, in the second case,
$N=T^lM'T^{-l}$ for some $l \in \mathbb Z$, hence $N \sim M$.

\end{pf}
\section{Further remarks}\label{s:furtherremarks}
\subsection{Determining Markoff matrices}\label{ss:determiningMM}
We saw that for a Markoff matrix $M=\left(%
\begin{array}{cc}
  a & b \\
  c & d \\
\end{array}%
\right)\in G$, $a+d=3|c|$. Nonetheless this is not a sufficient
condition, if we denote the conjugacy
class of $M$ in $\SLTwoZ$ by $cl(M)$, and let $M'=\left(%
\begin{array}{cc}
  a' & b' \\
  c' & d' \\
\end{array}%
\right)$ be the elements of $cl(M)$, then it follows from some
elementary hyperbolic geometry that
$$\min_{M'\in cl(M)} |c'|=|c|.$$
This extra information may allow us to obtain the Markoff
conjecture for more general $c$ using our methods, namely, even if
$M_k\in \SLTwoZ$, it may not be in $G$, or may not satisfy the
condition above, hence, may not be a Markoff matrix.

\subsection{Representatives of the equivalence classes of Markoff matrices}
\label{ss:reps} Representatives for the equivalence classes of the
Markoff matrices  can be enumerated  starting from $A$ and $B$ via
the Farey sequence in the following way: Associate the rational
number $0/1$ to the matrix $A$ and the number $1/1$ to the matrix
$B$. So using $W$ to denote the matrix corresponding to a word in
$A$ and $B$ we can write $W_0=A$, $W_1=B$. Define $W_{1/2}=AB$.
More generally if $p/q$ and $r/s$ are Farey neighbours with $0\le
p/q \le r/s<1$ then we define
$$W_{\frac{p+q}{r+s}}=W_{p/q}\cdot W_{r/s}$$
In other words, the word $W$ in $A$ and $B$ associated to the
rational $m/n$ is obtained by concatenating the words associated
to the previous two Farey numbers, the order of the concatenation
being determined by the ordering of the fractions. The
representatives for the equivalence classes of Markoff matrices
can be taken to be the matrices $M_q$ where $q \in {\mathbb Q}
\cap [0,1/2]$. Details are left to the reader.

\subsection{Relation to the gaps of McShane's and his
identity}\label{ss:Mcshane} There is an interesting relation
between our method and McShane's identity \cite{Mcshane}. Each
Markoff matrix (up to taking inverse and conjugation by $T$) gives
rise to an open interval/gap in ${\mathbb R}/{\mathbb Z}$ centered
at $a/c$, and of width $(3-\frac{\sqrt{9c^2-4}}{c})$, (where we
assume $c>0$). These gaps are all disjoint, their complement is a
Cantor set of measure 0, the fact that the sum is 1 is a special
case of McShane's identity (which works for all cusped hyperbolic
tori). It should be possible to exploit the disjointedness of
these gaps to provide stronger evidence for the Markoff
conjecture. The idea is that a shift by $k/c$ may result in a
matrix $M_k$ in $\SLTwoZ$ but the resulting gap may not be
disjoint from the other gaps already obtained.

\subsection{Other results}\label{ss:otherresults}
Our methods easily extend to give uniqueness for the case where
$c=2p^n$, and, using the conditions of \S \ref{ss:determiningMM}
and some slightly tedious calculations, the case $c=5p^n$. Baragar
\cite{Baragar} has also proven the Markoff conjecture for certain
other classes of Markoff numbers, his methods however are
different from ours, and our method does not apply to the classes
he considers.

{}
\end{document}